\documentclass[11pt]{article}
\usepackage{amssymb}
\newtheorem{theo}{Theorem}
\newtheorem{prop}{Proposition}

\newtheorem{coro}{Corollary}
\newtheorem{rema}{Remark}

\newcommand{\cqfd}
{%
\mbox{}%
\nolinebreak%
\hfill%
\rule{2mm}{2mm}%
\medbreak%
\par%
}
\date{\empty}
\title{A counterexample to the Hodge conjecture for\\
 K\"ahler
varieties}
\author{Claire Voisin\\ 
Institut de math\'ematiques de Jussieu, CNRS,UMR 
7586}
\begin{document}
\maketitle
\section{Introduction}\label{sec0} 
If $X$ is a smooth projective variety, it is in particular a K\"ahler variety,
and its cohomology groups carry the Hodge decomposition
$$ H^k(X,{\mathbb C})=\oplus_{p+q=k}H^{p,q}(X).$$
A class $\alpha\in H^{2p}(X,{\mathbb Q})$ is said to be a rational Hodge
class if its image in $H^{2p}(X,{\mathbb C})$ belongs to
$H^{p,p}(X)$. As is well known, the classes which are
Poincar\'e dual to irreducible algebraic subvarieties of codimension
$p$ of $X$ are degree $2p$
Hodge classes. The Hodge conjecture asserts that any rational 
Hodge class
is a combination with rational coefficients of such classes.

\par In the case of a general compact K\"ahler variety $X$, 
the conjecture
 above, where the algebraic subvarieties are replaced 
 with closed analytic
subsets, is known to be false (cf \cite{zu}). The simplest 
example is provided by a  complex
torus
endowed with a holomorphic line bundle of indefinite curvature.
 If the
torus is chosen general enough, it will not
contain any analytic hypersurface, while the first Chern class
of the line bundle will provide a Hodge class of degree $2$.

\par In fact, another general method to construct
Hodge classes is to consider Chern classes of holomorphic
vector bundles. In the projective case, the set of classes
generated this way is the same as the set generated by classes of
subvarieties. To see this, one looks at a still more general
set of classes, which is the set generated by the Chern classes
of coherent sheaves on $X$. Since any coherent sheaf has a finite
 resolution by locally free sheaves, one does not get more classes
 than with locally free sheaves. On the other hand, this later
 set obviously contains the classes of subvarieties
 (one computes for this the top Chern class of
 ${\mathcal I}_Z$ for $Z$ irreducible of codimension $p$, 
 and one shows that it is proportional
to the class of $Z$). Finally, to see that
the classes of coherent sheaves can be generated by classes
of subvarieties, one puts a filtration on any
coherent sheaf, whose associated graded consists of 
rank $1$ sheaves supported on subvarieties, which makes the
result easy.
\par In the general K\"ahler case, none of these equalities holds.
The only
obvious result is that the space generated by the Chern classes of
analytic coherent sheaves contains  both
the classes which are
Poincar\'e dual to irreducible  closed analytic subspaces
and the Chern classes of holomorphic vector bundles
(or locally free analytic coherent sheaves).
The example
above shows that a Hodge class of degree $2$ may be the Chern class 
of a holomorphic line bundle, even if $X$
 does not contain any complex 
analytic subset. On the other hand, 
it may be the case that coherent sheaves do not admit
a resolution by locally free sheaves, (although
it is true in dimension $2$ \cite{sch}), and that
more generally $X$ does not carry enough vector bundles to generate
the Hodge classes of subvarieties or coherent sheaves
(see the appendix for such examples).
 Hence the set generated
by the Chern classes of analytic coherent
sheaves is actually larger than the two others.

\par Notice that using the
Grothendieck-Riemann-Roch
formula  (cf \cite{BS},
extended in \cite{OTT}
to the complex analytic case),
 we can give the following alternative description
of this set : it is generated by the classes
$$\phi_*c_i(E)$$
where $\phi:Y\rightarrow X$
is a morphism from another compact K\"ahler manifold,
$E$ is a holomorphic vector bundle on $Y$, and 
$i$ is any integer. Another fact which follows
from iterated applications of the Whitney formula, is that
the set which is additively generated by
the Chern classes of coherent sheaves is equal to the set which is
generated as a subring of the cohomology ring
by
the Chern classes of coherent sheaves. Thus this set is as big
and stable as possible.

\par 
 Now, since we do not know other
 geometric ways of constructing Hodge classes, the following seems to
 be a natural extension 
of the Hodge conjecture to K\"ahler varieties.

\par{\it Are the rational Hodge classes of
a compact K\"ahler variety $X$  generated over ${\mathbb Q}$ by Chern
classes of analytic coherent  sheaves on $X$?}

\par Our goal in this paper is to give a negative answer
to this question. We show the
following theorem

\begin{theo} 
\label{intro}There exists a $4$-dimensional complex torus $X$ 
which possesses 
a non trivial Hodge class of degree $4$, such that any 
analytic coherent sheaf ${\mathcal F}$ on $X$
satisfies $c_2({\mathcal F})=0$.
\end{theo}
In the appendix, we also give a few geometric
consequences of a result of Bando and Siu \cite{Siu}, extending Uhlenbeck-Yau's
theorem. We show in particular 
that for a general  compact K\"ahler variety
$X$, the analytic coherent sheaves
on $X$ do not admit finite resolutions by
locally free coherent sheaves. This answers a question asked to us
by L. Illusie.

{\bf Acknowledgements.} I would like to thank Joseph Le Potier and
Andrei Teleman for helpful discussions on this paper.  The references
\cite{bapo} and \cite {te} have been a starting point for this work.
I also thank Pierre Deligne and Luc Illusie for their interest and their questions.
\section{A criterion for the vanishing of  Chern
 classes of coherent sheaves \label{sec2}}
We consider in this section a compact K\"ahler variety $X$ 
of dimension $n\geq3$ satisfying the following assumptions
\begin{enumerate}
\item \label{a} The N\'eron-Severi group $NS(X)$ generated by the first Chern
classes of holomorphic line bundles on $X$ is equal to $0$.
\item \label{b} $X$ does not contain any proper closed analytic subset of
positive dimension.
\item\label{c} For some K\"ahler class $[\omega]\in H^2(X,{\mathbb R})\cap
H^{1,1}(X)$, the set of Hodge classes $Hdg^4(X,{\mathbb Q})$ is 
perpendicular to $[\omega]^{n-2}$ for the intersection pairing
$$H^4(X,{\mathbb R})\otimes H^{2n-4}(X,{\mathbb R})
\rightarrow{\mathbb R}.$$
\end{enumerate}
Our aim is to show the following
\begin{prop}\label{prop1} If $X$ is as above, any analytic
 coherent sheaf ${\mathcal F}$ on
$X$ satisfies $c_2({\mathcal F})=0$.
\end{prop}
{\bf Proof.} As in \cite{bapo}, the proof is by induction on the rank $k$ of 
${\mathcal F}$. We note that because $dim\,X\geq3$, torsion sheaves
supported on points on $X$ have trivial $c_1$ and $c_2$.
On the other hand
assumption \ref{b} implies that torsion sheaves are supported
 on points. This deals with the case
 where $rk\,{\mathcal F}=0$. Furthermore this shows 
 that we can restrict
 to torsion free coherent sheaves.

\par If ${\mathcal F}$
contains a non trivial subsheaf ${\mathcal G}$ of rank $<k$,
we have the exact sequence
$$0\rightarrow {\mathcal G}\rightarrow{\mathcal F}
\rightarrow{\mathcal F}/{\mathcal G}\rightarrow 0,$$
with $rk\,{\mathcal G}<k$ and 
$rk\,{\mathcal F}/{\mathcal G}<k$. Assumption \ref{a}
and the induction hypothesis then give
$$c_1({\mathcal G})=c_2({\mathcal G})=0,\,
c_1({\mathcal F}/{\mathcal G})=c_2({\mathcal F}/{\mathcal G})=0.$$
Therefore Whitney formula implies $c_2({\mathcal F})=0$.

\par We are now reduced to study the case where ${\mathcal F}$ 
does not contain
any non trivial proper subsheaf of
smaller rank. By assumption 
\ref{b}, ${\mathcal F}$ is locally free away from
 finitely many points
$\{x_1,\ldots,x_N\}$ of $X$.
One shows now
that there exists
a variety $$
\tau:\tilde X\rightarrow X,$$ which is obtained from $X$ by finitely 
many successive
blow-ups with smooth centers (and in particular is also 
K\"ahler), so that $\tau$ restricts to an isomorphism
$$\tilde X-\cup_i E_i\cong X-\{x_1,\ldots,x_N\},$$
where $E_i:=\tau^{-1}(x_i)$
and such that there exists a locally free
sheaf $\tilde{\mathcal F}$ on $\tilde X$ which is isomorphic to
${\mathcal F}$ on the open set $\tilde X-\cup_i E_i$ :
indeed, the problem is local near each $x_i$.
Choosing a finite presentation
$${\mathcal O}_X^l\rightarrow {\mathcal O}_X^r
\rightarrow {\mathcal F}\rightarrow 0$$
of ${\mathcal F}$ near $x_i$, we get a morphism
to the Grassmannian of $l$-dimensional
subspaces of ${\mathbb C}^r$, which is well defined away
from $x_i$, since ${\mathcal F}$ is free
away from $x_i$.
This morphism is easily seen to be meromorphic. Hence by Hironaka
desingularization theorem \cite{hi}, this morphism can be
extended after finitely many blowups. Then
the pull-back of the tautological quotient
bundle on the Grassmannian will provide the desired extension.

\par
Note that because ${\mathcal F}$ does not contain any
non zero
subsheaf of smaller rank, the same is true of $\tilde{\mathcal F}$, 
which
implies that
$\tilde{\mathcal F}$ is $h_\lambda$-stable for any
K\"ahler metric $h_\lambda$ on
$\tilde X$. The theorem of existence of Hermitian-Yang-Mills 
connections
(\cite{UY}) then provides $\tilde{\mathcal F}$
with  a  Hermitian-Einstein metric
$k_\lambda$ for any K\"ahler metric $h_\lambda$
on $\tilde X$. This means 
 that the curvature $R_\lambda\in \Gamma({\mathcal
H}om\,(\tilde{\mathcal F},
\tilde{\mathcal F})\otimes\Omega_{\tilde X,{\mathbb R}}^2)$
of the metric connection on
$\tilde{\mathcal F}$ associated to  $k_\lambda$
 is the sum of
a diagonal matrix with all coefficients equal to
$\mu_\lambda\omega_\lambda$
and of a matrix $R_\lambda^0$
whose coefficients are $(1,1)$-forms
anihilated by the $\Lambda$ operator
relative to the metric $h_\lambda$. (The connection
is then said to be Hermitian-Yang-Mills.)
Here $\omega_\lambda$ is the K\"ahler form of the metric
$h_\lambda$ and $\mu_\lambda$ is a constant coefficient, equal to
\begin{eqnarray}
\label{eqpre}
2i\pi\frac{c_1(\tilde{\mathcal
F})[\omega_\lambda]^{n-1}}{k[\omega_\lambda]^{n}},
\end{eqnarray}
where $[\omega_\lambda]\in H^2(\tilde X,{\mathbb R})$ denotes the
de Rham class of the form
$\omega_\lambda$.
\par  Let us denote by $\eta_\lambda^0$ the closed $4$-form
$$\eta_\lambda^0=tr\,(\frac{R_\lambda^0}{2i\pi})^2.$$
We assume chosen small neighbourhoods
$V_i$ of $x_i$ in $X$, and forms $\omega_i$
on $\tau^{-1}(V_i)$ which vanish near the boundary of
$\tau^{-1}(V_i)$, and restrict to a K\"ahler form
on $\tau^{-1}(x_i)$. Then we will choose
\begin{eqnarray}
\label{jennifer}
\omega_\lambda=\tau^*\omega+\lambda(\sum_i\omega_i)
\end{eqnarray}
which is easily seen to be K\"ahler for
sufficiently small $\lambda$. We shall now prove
that $R_\lambda$ tends to $0$ with $\lambda$ in the
$L^2$-sense away from the $V_i$'s.
The argument is an extension
of L\"ubke's inequality \cite{Lu} which proves that  a
Hermitian-Yang-Mills connection on a vector bundle $E$ with
$c_1(E)[\omega]^{n-1}=c_1(E)^2[\omega]^{n-2}=c_2(E)[\omega]^{n-2}=0$,
 where
$[\omega]$ is the class of the K\"ahler form on the basis,
is in fact
 flat.
We first claim :
\begin{prop}\label{jeudi}
For any  differential $2n-4$-form $\alpha$ on
$X$, the integral
$$\int_{X-\cup_iV_i}\eta_\lambda^0\wedge\alpha$$
tends to $0$ with $\lambda$.
\end{prop}
Before proving this proposition, we show how it implies
that $c_2({\mathcal F} )=0$, thus completing the
induction step.

\par Poincar\'e duality will provide an isomorphism
$$H^4(X-\cup_iV_i,{\mathbb R})\cong
H^{2n-4}(X,\cup_iV_i,{\mathbb R})^*,$$
which is realized by integrating  closed
$4$-forms defined over $X-\cup_iV_i$ against
$2n-4$-forms vanishing on the $V_i$'s.
Next because $dim\,X\geq 3$ we have isomorphisms
$$H^4(X,{\mathbb R})\cong H^4(X-\cup_iV_i,{\mathbb R}),$$
$$H^{2n-4}(X,\cup_iV_i,{\mathbb R})\cong
H^{2n-4}(X,{\mathbb R})$$
which are compatible with Poincar\'e duality.
Now the cohomology class of the closed $4$-form $\eta_\lambda^0$
is easily computed to be 
\begin{eqnarray}
\label{marnounou}
[\eta_\lambda^0]=-2c_2(\tilde{\mathcal F})+c_1(\tilde{\mathcal F})^2
-\frac{\mu_\lambda}{i\pi}[\omega_\lambda]\cup 
c_1(\tilde{\mathcal F})+
k(\frac{\mu_\lambda}{2i\pi})^2[\omega_\lambda]^2.\end{eqnarray}
Hence its restriction to
$\tilde X-\cup_i\tau^{-1}(V_i)=X-\cup_iV_i$ is equal
to $-2c_2({\mathcal F})+
k(\frac{\mu_\lambda}{2i\pi})^2[\omega]^2$, since
$c_1({\mathcal F})=0$ and $\omega_\lambda$ restricts to
$\omega$ on $X-\cup_iV_i$.

In order to show that $c_2({\mathcal F})=0$,
it then suffices by the above to show that for any
closed $2n-4$-form
$\alpha$
on $X$, vanishing on
$\cup_i V_i$, we have
$$\int_X(\eta_\lambda^0-k(\frac{\mu_\lambda}{2i\pi})^2\omega^2)
\wedge\alpha=0.$$
But this integral is independent of
$\lambda$ and so it suffices to show
that 
\begin{eqnarray}
\label{nouveeq}\lim_{\lambda\rightarrow 0}
\int_X(\eta_\lambda^0
-k(\frac{\mu_\lambda}{2i\pi})^2\omega^2)\wedge\alpha=0.
\end{eqnarray}
Now we claim that
 \begin{eqnarray}
 \label{nomlun}
 \lim_{\lambda\rightarrow 0}\mu_\lambda=0.
 \end{eqnarray}
Indeed this follows from formula
(\ref{eqpre}), and from the fact that the class
$c_1(\tilde {\mathcal F})$ restricts to $0$ on
$\tilde X-\cup_i\tau^{-1}(x_i)$ because
$NS(X)=0$.  Then the intersection pairing
$<c_1(\tilde {\mathcal F}),\tau^*[\omega]^{n-1}>_{\tilde X}$
is equal to $0$,
and we conclude using the fact that
$\lim_{\lambda\rightarrow 0}\omega_\lambda=\tau^*\omega$. 

Then (\ref{nouveeq}) follows from (\ref{nomlun})
and Proposition \ref{jeudi}.

\par

We now go to the proof of proposition
\ref{jeudi}. 
We  first claim that
\begin{eqnarray}
\label{marnou}
\lim_{\lambda\rightarrow 0}
\int_{\tilde X}\eta_\lambda^0\wedge
\omega_\lambda^{n-2}=0.\end{eqnarray}
Indeed, we know that the space
$Hdg^4(X) $ is perpendicular to
$[\omega]^{n-2}$ for the intersection pairing.
Hence we have
$$<c_2(\tilde {\mathcal F}),\tau^*[\omega]^{n-2}>_{\tilde X}=0.$$
On the other hand this is equal to
$$\lim_{\lambda\rightarrow 0}
<c_2(\tilde {\mathcal F}),[\omega_\lambda]^{n-2}>_{\tilde X}$$
since $\lim_{\lambda\rightarrow 0}\omega_\lambda=\tau^*\omega$.
Exactly by the same argument, we show that
$$\lim_{\lambda\rightarrow 0}<
c_1^2(\tilde {\mathcal F}),[\omega_\lambda]^{n-2}>_{\tilde X}=0.$$
Then the result follows from formula
(\ref{marnounou}) and from (\ref{nomlun}).
\cqfd
Next we recall  that the endomorphism
$R_\lambda^0$  of
$\tilde{\mathcal F}$, with forms coefficients is
anti-self-adjoint with respect to the metric
$k_\lambda$. This follows from the fact that 
$R_\lambda$ is the curvature of the
metric  connection with respect to $k_\lambda$.
In a local orthonormal basis of
$\tilde{\mathcal F}$, this will be translated into the fact
that $R_\lambda^0$ is represented by a matrix,
whose coefficients are differential $2$-forms,
 which satisfies
$$\overline{^tR_\lambda^0}=-R_\lambda^0.$$
The second crucial property of
$R_\lambda^0$ is the fact that its coefficients are
primitive differential $(1,1)$-forms on
$\tilde X$, with respect to the metric
$h_\lambda$. It is well known that
this implies the following equality 
$$*_\lambda\overline{\gamma}= -
\overline{\gamma}\wedge\frac{\omega_\lambda^{n-2}}{(n-2)!},$$
where $*_\lambda$ is the Hodge $*$-operator
for the metric $h_\lambda$. Since
$h_\lambda$ restricts to $h$ on
$X-\cup_iV_i$, these forms satisfy as well
\begin{eqnarray}
\label{eqvencre}
*\overline{\gamma}=-\overline{\gamma}\wedge
\frac{\omega^{n-2}}{(n-2)!}
\end{eqnarray}
on $X-\cup_iV_i$.

\par
Now let $\alpha$ be a differential $2n-4$-form on $X$.
Then it follows from (\ref{eqvencre}) that
there exists a positive constant $c_\alpha$ such that for any
primitive $(1,1)$-form $\gamma$ on $X$, we have the 
following pointwise inequality of pseudo-volume forms
on $X$:
\begin{eqnarray}
\label{eqeq}\mid\gamma\wedge\overline{\gamma}\wedge\alpha
\mid\leq c_\alpha\gamma\wedge*\overline{\gamma}
=-c_\alpha\gamma\wedge\overline{\gamma}\wedge\frac
{\omega^{n-2}}{(n-2)!}.
\end{eqnarray}
Working locally in
a orthonormal basis of $\tilde{\mathcal F}$ and
using the fact that the matrix $R_\lambda^0$ is
anti-self-adjoint and with primitive coefficients
of $(1,1)$-type, we now
get
the pointwise inequality of pseudo-volume
 forms on $X-\cup_iV_i$
$$\mid tr\,(R_\lambda^0)^2\wedge\alpha\mid
\leq c_\alpha tr\,(R_\lambda^0)^2\wedge\frac{\omega^{n-2}}{(n-2)!}.$$
Therefore we get
the inequality
$$\mid\int_{X-\cup_iV_i}tr\,(R_\lambda^0)^2\wedge\alpha\mid
\leq c_\alpha\int_{X-\cup_iV_i}
tr\,(R_\lambda^0)^2\wedge\frac{\omega^{n-2}}{(n-2)!}.$$
But by (\ref{marnou}), and because
$\eta_\lambda^0=Tr\,(\frac{R_\lambda^0}{2i\pi})^2$,
 we know
that 
$$\lim_{\lambda\rightarrow 0}
\int_{\tilde X}Tr\,(R_\lambda^0)^2\wedge
\omega_\lambda^{n-2}=0.$$
Because the integrand is positive
and $\omega=\omega_\lambda$ on $X-\cup_iV_i$, this implies that
$$\lim_{\lambda\rightarrow 0}
\int_{X-\cup_iV_i}Tr\,(R_\lambda^0)^2\wedge
\omega^{n-2}=0.$$
Hence 
$$\lim_{\lambda\rightarrow 0}\int_{X-\cup_iV_i}
tr\,(R_\lambda^0)^2\wedge\alpha=0=
\lim_{\lambda\rightarrow 0}\int_{X-\cup_iV_i}
\eta_\lambda^0\wedge\alpha.$$
Proposition \ref{jeudi} is proven.
\cqfd
\begin{rema} A. Teleman mentioned to me the possibility of using
the result of \cite{Siu} (see the appendix) to give a shorter proof of
the equality $c_2({\mathcal F})=0$ for stable ${\mathcal F}$. In this paper, the
results of Uhlenbeck and Yau are extended to reflexive coherent stable
sheaves, and L\"ubke's inequality, together
with the fact that equality implies
projective flatness, is proven.

Since in our case we have a much stronger assumption than stability,
namely stability of any desingularization of ${\mathcal F}$
with respect to any K\"ahler metric, it seemed however
appropriate to avoid the reference to the technically hard
result of \cite{Siu} and to content ourselves with
an argument which appeals only to \cite{UY} and elementary
computations.
\end{rema}

\section{Constructing an example}
Our example will be of Weil type \cite{W}. The Hodge classes
described below have been constructed by Weil in the
case of algebraic tori, as a potential counterexample
to the Hodge conjecture for algebraic varieties.
In the case of a general complex torus, the construction is still
simpler.
These complex tori have been also considered in
\cite{zu} by Zucker , who proves there some of the 
results stated below. (I thank
 P. Deligne and C. Peters for pointing
out this reference to me.)

\par We start with a ${\mathbb Z}[I]$-action
on $\Gamma:={\mathbb Z}^8$, where $I^2=-1$,
which makes
$$\Gamma_{\mathbb Q}:=\Gamma\otimes{\mathbb Q}$$
into a $K$-vector space, where
$K$ is the quadratic field ${\mathbb Q}[I]$.

Let 
$$\Gamma_{\mathbb C}=\Gamma\otimes{\mathbb C}={\mathbb C}^4_i\oplus
{\mathbb C}^4_{-i}$$
be the associated decomposition into eigenspaces for $I$.
 A four dimensional complex torus $X$ with underlying lattice
  $\Gamma$ and inheriting the $I$-action is determined by
  a $4$ dimensional complex subspace
  $W$ of $\Gamma_{\mathbb C}$, which has to be the direct
  sum 
  $$W=W_i\oplus W_{-i}$$
  of its intersections with ${\mathbb C}^4_i$
  and ${\mathbb C}^4_{-i}$.  It has furthermore to satisfy the condition
  that 
  \begin{eqnarray}
  \label{cond}W\cap\Gamma_{\mathbb R}=\{0\}.
  \end{eqnarray}
  Given $W$, $X$ is given by the formula
  $$X=\Gamma_{\mathbb C}/(W\oplus\Gamma).$$
  
  We will choose $W$ so that 
  $$dim\,W_i=dim\,W_{-i}=2.$$
  Then $W$, hence $X$ is determined by the choice of the 
  $2$-dimensional subspaces
  $$W_i\subset {\mathbb C}^4_i,\,W_{-i}\subset{\mathbb C}^4_{-i},$$
  which have to be general enough so that
  the condition (\ref{cond})
  is satisfied.
  
  We have isomorphisms
  \begin{eqnarray}
  \label{soir}
  H^4(X,{\mathbb Q})\cong H_4(X,{\mathbb Q})\cong
  \bigwedge^4\Gamma_{\mathbb Q}.
   \end{eqnarray}
  Consider the subspace
  $$\bigwedge_K^4\Gamma_{\mathbb Q}\subset
  \bigwedge^4\Gamma_{\mathbb Q}.$$
  Since $\Gamma_{\mathbb Q}$ is a $4$-dimensional
  $K$-vector space, $\bigwedge_K^4\Gamma_{\mathbb Q}$
   is a one dimensional $K$-vector space, and its image
   is a $2$ dimensional ${\mathbb Q}$-vector space.
   The claim is that $\bigwedge_K^4\Gamma_{\mathbb Q}$
   is made of Hodge classes, that is, is contained in the
   subspace $H^{2,2}(X)$ for the Hodge decomposition.
   Notice that under the isomorphisms
   (\ref{soir}), tensorized by
   ${\mathbb C}$, 
   $H^{2,2}(X)$ identifies to the image
   of
   $$\bigwedge^2W\otimes
   \bigwedge^2\overline{W}$$
   in $\bigwedge^4\Gamma_{\mathbb C}.$

   To prove this claim, note that we have the decomposition
   $$\Gamma_K:=\Gamma_{\mathbb Q}\otimes K=\Gamma_{K,i}
   \oplus\Gamma_{K,-i}$$
   into eigenspaces for the $I$ action.
   Then
   $\bigwedge^4_K\Gamma_{\mathbb Q}\subset
   \bigwedge^4\Gamma_{\mathbb Q}$
   is defined as  the image
   of $\bigwedge_K^4\Gamma_{K,i}\subset
   \bigwedge_K^4\Gamma_K$ via the trace map
   $$\bigwedge_K^4\Gamma_K=
   \bigwedge_{\mathbb Q}^4\Gamma_{\mathbb Q}\otimes K
   \rightarrow
    \bigwedge^4\Gamma_{\mathbb Q}.$$
    
    Now we have the inclusion
    $$\Gamma_K\subset\Gamma_{\mathbb C}$$
    and the equality
    $$\Gamma_{K,i}=\Gamma_K\cap {\mathbb C}^4_i.$$
    The space
   $ \Gamma_{K,i}$ is a $4$ dimensional $K$-vector space
   which generates over ${\mathbb R}$ the space
    ${\mathbb C}^4_i$. It follows that the image
   of $\bigwedge_K^4\Gamma_{K,i}$
   in $\bigwedge^4\Gamma_{\mathbb C}$ generates
   over ${\mathbb C}$
   the line
   $\bigwedge^4{\mathbb C}^4_i$.
   
   \par But we know that ${\mathbb C}^4_i$ is the direct sum
   of the two spaces
   $W_i$ and $\overline{W_{-i}}$ which are 
   $2$-dimensional.
   Hence 
   $$\bigwedge^4{\mathbb C}^4_i=\bigwedge^2W_i\otimes
   \bigwedge^2\overline{W_{-i}}$$
   is contained in
   $\bigwedge^2W\otimes\bigwedge^2\overline{W}$,
   that is in $H^{2,2}(X)$.
   \cqfd
   To conclude the construction of an example satisfying the conclusion of
   proposition \ref{prop1},
   and hence the proof of theorem
   \ref{intro}, it remains now only to prove that a general $X$
   as above satisfies the assumptions stated at the
   beginning of section \ref{sec2}. Since $X$
   is a complex torus, assumption \ref{b}
   will be a consequence of assumption \ref{a} and of the fact
   that $X$ is simple. Indeed it is known that
   if $Y\subset X$ is a proper positive dimensional
   subvariety of a simple complex torus, then
   $Y$ has positive canonical bundle. But $X$ being simple,
   $Y$ must generate $X$ as a group, and then
   $X$ must be algebraic, contradicting the
   fact that $NS(X)=0$.
   
   Next we show that the Hodge classes in
    $\bigwedge^4_K\Gamma_{\mathbb Q}$ constructed
    above are perpendicular to $[\omega]^2$
    for a K\"ahler class
    $[\omega]\in H^{1,1}(X)$. To see this, note
     that with the notations
    as above, these classes ly in
    $\bigwedge^2W_i\otimes\bigwedge^2\overline{W_{-i}}$,
    with
    $$W_i\subset W,\,\overline{W_{-i}}\subset\overline{W}.$$
    The space $H^{1,0}(X)$ identifies to
    $\overline{W}^*$ and accordingly the space
    $H^{1,1}(X)$ identifies to
    $\overline{W}^*\otimes W^*$. For $[\omega]\in
    \overline{W}^*\otimes W^*$, the pairing
    $<[\omega]^2,H^{2,2}(X)>$, restricted to
    $\bigwedge^2W_i\otimes\bigwedge^2\overline{W_{-i}}$,
     is obtained by squaring $[\omega]$
     to get an element of
     $$\bigwedge^2\overline{W}^*\otimes\bigwedge^2 W^*
     \cong\bigwedge^2 W^*\otimes\bigwedge^2\overline{W}^*$$
     and by projecting to
     $\bigwedge^2W_i^*\otimes\bigwedge^2\overline{W_{-i}}^*$.
     
     \par Now choose 
     $$[\omega]\in \overline{W_i}^*\otimes W_i^*
     \oplus\overline{W_{-i}}^*\otimes W_{-i}^*.$$
     Since
     $\overline{W}^*=\overline{W_i}^*\oplus
     \overline{W_{-i}}^*$, we can find a 
     K\"ahler class $[\omega]$ in this space.
     On the other hand, we see that $[\omega]^2$
     belongs to the space
     $$\bigwedge^2\overline{W_i}^*\otimes\bigwedge^2W_i^*
     \oplus \overline{W_i}^*\otimes
     \overline{W_{-i}}^*\otimes W_i^*\otimes W_{-i}^*
     \oplus
     \bigwedge^2\overline{W_{-i}}^*\otimes\bigwedge^2W_{-i}^*.$$
     Hence its projection (after switching the factors in the tensor
     product) to
     $\bigwedge^2W_i^*\otimes\bigwedge^2\overline{W_{-i}}^*$
     is equal to $0$.
     \cqfd
     
     In conclusion, the assumptions at the beginning
     of section \ref{sec2}
     will be a consequence of the following
     facts : 
     \begin{prop} \label{profi}
      For a general $X$ as above, we have :
      \begin{enumerate}
     \item\label{asam} $NS(X)=0$.
     \item\label{bsam} $X$ is simple.
     \item\label{csam} The space $Hdg^4(X)$ is equal
     to 
     the space 
     $$\bigwedge^4_K\Gamma_{\mathbb Q}
     \subset \bigwedge^4\Gamma_{\mathbb Q}=H_4(X,{\mathbb Q})\cong H^4(X,{\mathbb Q}).$$
     \end{enumerate}
     \end{prop}
     {\bf Proof.}
      The analogues of these statements have been 
     proven in the algebraic case in \cite{W} (see also
     \cite{VG}).The result is that
     for a general abelian $4$-fold of
     Weil type, the N\'eron-Severi group
     is of rank $1$, generated by a class $\omega$, and 
     the space $Hdg^4(X)$ is of rank $3$, generated over
     ${\mathbb Q}$ by
     the space 
     $\bigwedge^4_K\Gamma_{\mathbb Q}
     $ and by the class
     $\omega^2$.
     Furthermore property \ref{bsam} is true
     for the generic abelian variety $X$ of Weil type.
     
     \par Property \ref{bsam} for the general
     complex torus of Weil type follows
     immediately, since this is a property satisfied
     away from the countable union
     of closed analytic subsets of the
     moduli space of complex tori of Weil type.
     
     \par As for properties \ref{asam} and \ref{csam}, we prove them
     by an infinitesimal argument, starting
     from an abelian 4-fold of Weil type
     $X$ satisfying
     the properties stated above.
     Assume  we can show that for some first order
     deformation
     $u\in H^1(T_X)$, tangent to the moduli
     space of complex tori
      of Weil type (which is smooth), we have
      $$int(u)(\omega)\not=0\,{\rm in}\, H^{2}({\mathcal O}_X),$$
      where the interior product here is composed of the cup-product
      $$H^1(T_X)\otimes H^1(\Omega_X)\rightarrow 
      H^2(T_X\otimes\Omega_X)$$
      and of the map induced by  the contraction
      $$ H^2(T_X\otimes\Omega_X)\rightarrow H^2({\mathcal O}_X).$$
      Then from the general theory of  Hodge loci
      (\cite{GrH}), it will
      follow that for a general complex torus of Weil
      type, we have $NS(X)=0$.
      Furthermore,  it will also follow
      that
       $$int(u)(\omega^2)\not=0\,{\rm in}\, H^{3}(\Omega_X),$$
       because it is equal to
       $2\omega\cup int(u)(\omega)$ and the
       cup product with $\omega$
       from
       $H^2({\mathcal O}_X)$ to
       $H^3(\Omega_X)$ is injective
       because $\omega$ is a K\"ahler class
      on $X$. But as before this will imply
      by the theory of Hodge loci that for a general complex
      torus of Weil type we have
      $rk\,Hdg^4(X)=2$ so that
      $Hdg^4(X)=\bigwedge^4_K\Gamma_{\mathbb Q}$.
      
      \par Hence it remains only to find
      such $u$, which is equivalent to prove that
      if for any $u$ tangent to the moduli space
      of complex tori of Weil type,
      $int(u)(\omega)=0$
      in $H^2({\mathcal O}_X)$, then $\omega=0$ in $H^1(\Omega_X)$.
      Here the notations are as in the beginning of this section.
      The tangent space to the deformations
      of the complex torus $X$ identifies
      to 
      \begin{eqnarray}
      \label{Tgt}
       Hom\,(W,\Gamma_{\mathbb C}/W)=Hom\,(W,\overline{W})=W^*\otimes
      \overline{W}.
      \end{eqnarray}
      The tangent space to the deformations of $X$ as a complex torus
      of Weil type is
      then the subspace
      $$Hom\,(W_i,{\mathbb C}^4_i/W_i)\oplus
      Hom\,(W_{-i},{\mathbb C}^4_{-i}/W_{-i})$$
      $$=W_i^*\otimes\overline{W_{-i}}\oplus
      W_{-i}^*\otimes\overline{W_i}.$$
      Via the identification
      (\ref{Tgt}), the interior
      product
      $$H^1(X,T_X)\otimes H^1(X,\Omega_X)\rightarrow 
      H^2(X,{\mathcal O}_X)=\bigwedge^2{W}^*$$
      identifies to the
      contraction followed by the wedge product
      $$W^*\otimes
      \overline{W}\otimes\overline{W}^*\otimes W^*
      \rightarrow \bigwedge^2W^*.$$
      We now write
      $$\omega=\omega_1+\omega_2+\omega_3+\omega_4,$$
      where
      $$\omega_1\in\overline{W_i}^*\otimes W_i^*,\,
      \omega_2\in\overline{W_i}^*\otimes W_{-i}^*,$$
      $$\omega_3\in\overline{W_{-i}}^*\otimes W_i^*,\,
      \omega_4\in\overline{W_{-i}}^*\otimes W_{-i}^*.$$
      Then clearly
      for $u_1\in W_i^*\otimes\overline{W_{-i}}$ we have
      $$int(u_1)(\omega_1)=int(u_1)(\omega_2)=0,$$
      $$ int(u_1)(\omega_3)\in\bigwedge^2W_i^*,\,
       int(u_1)(\omega_4)\in W_i^*\otimes W_{-i}^*.$$
       Similarly, for $u_2\in W_{-i}^*\otimes\overline{W_i}$
       we have
       $$int(u_2)(\omega_3)=int(u_1)(\omega_4)=0,$$
      $$ int(u_2)(\omega_2)\in\bigwedge^2W_{-i}^*,\,
       int(u_2)(\omega_1)\in W_i^*\otimes W_{-i}^*.$$
      
     The condition 
     $$int(u_1)(\omega)=0=
     int(u_2)(\omega)=0$$ for any
     $u_1,\,u_2$ then implies that
     $$int(u_1)(\omega_3)=0\,{\rm in}\, \bigwedge^2W_i^*,\,
     int(u_1)( \omega_4)=0\,
     {\rm in}\, W_i^*\otimes W_{-i}^*,$$
     $$int(u_2)(\omega_1)=0\,{\rm in}\,W_i^*\otimes W_{-i}^*
     ,\,
     int(u_2)(\omega_2)=0\,{\rm in}\,\bigwedge^2W_{-i}^*
     $$
     for any $u_1,\,u_2$.
     But it is obvious that it implies
     $\omega_1=\omega_2=\omega_3=\omega_4=0.$
     
     \cqfd

     Hence Proposition \ref{profi} is proven,
     which together with Proposition
     \ref{prop1} completes the proof
     of  Theorem \ref{intro}.
     \cqfd 
     
   \section{Appendix \label{appendix}}

Our goal in this appendix is to give a few geometric consequences
of the following theorem due to Bando and Siu
\cite{Siu} :
\begin{theo}
\label{siu} Let $X$ be a compact K\"ahler variety, 
endowed with a K\"ahler metric $h$ and let
${\mathcal F}$ be a reflexive  $h$-stable
sheaf on $X$. Then there
exists a Hermite-Einstein metric 
on ${\mathcal F}$ relative to $h$. Furthermore,
if we have  
 $<c_2({\mathcal F}),[\omega]^{n-2}>_X=0$
and  $c_1({\mathcal F})=0$, ${\mathcal F}$ is locally free
and
the associated metric connection is
flat.
\end{theo}
Here $[\omega]$ is the K\"ahler class of the
metric $h$.
\begin{rema} Once we know that the metric connection is flat
away from the singular locus $Z$ of ${\mathcal F}$, the fact that
${\mathcal F}$ is locally free is immediate.
Indeed the flat connection is associated to
  a local system on
  $X-Z$. But since $codim\,Z\geq2$, this local system extends
  to $X$. Hence there exists
  a holomorphic vector bundle $E$ on $X$, which admits
  a unitary flat connection and
   is isomorphic to
  ${\mathcal F}$ away from $Z$. 
  But since ${\mathcal F}$ is reflexive, the isomorphism
  $E\cong {\mathcal F}$ on $X-Z$ extends to $X$.
  \end{rema}

We assume now that $X$ is compact K\"ahler
and satisfies the condition
that the group $Hdg^2(X)$ of rational Hodge classes of degree $2$ vanishes,
and that the group
$Hdg^4(X)$ is perpendicular for the intersection pairing
to $[\omega]^{n-2}$ for some K\"ahler class
$\omega$ on $X$. Under these assumptions, $X$ does not contain any proper
analytic subset
of codimension less or equal to 2, and any
coherent sheaf ${\mathcal F}$ satisfies
the conditions
$$c_1({\mathcal F})=0,\,<c_2({\mathcal F}),[\omega]^{n-2}>_X=0.$$
We now prove :
\begin{prop}
\label{main} If $X$ is as above, for any 
torsion free coherent sheaf ${\mathcal F}$
on $X$, there exists a holomorphic  vector bundle
$E$ on $X$, whose all 
rational Chern classes $c_i(E),\,i>0$ vanish, and an exact sequence
$$0\rightarrow {\mathcal F}\rightarrow E\rightarrow {\mathcal T}\rightarrow 0$$
where ${\mathcal T}$ is a torsion sheaf on $X$.
\end{prop}

Before we prove the proposition, we state the following corollaries.
\begin{coro} \label{above} If $E$  is a holomorphic vector bundle on
$X$, then all rational Chern classes $c_i(E),\,i>0$ vanish.
\end{coro}
Indeed, we know that there
exists an inclusion
$$E\hookrightarrow E',$$
where $E'$ is a vector bundle of the same rank
as $E$ and satisfies the property that
 all rational Chern classes $c_i(E'),\,i>0$ vanish. Now since
 $Hdg^2(X)=0$, $X$ does not contain any 
 hypersurface, and it follows that the inclusion 
 above is an isomorphism.
 \cqfd
 \begin{coro} If $X$ is as above and
 $Z\subset X$ is a non-empty proper analytic subset, 
 the ideal sheaf
 ${\mathcal I}_Z$ does not admit a finite free resolution.
 \end{coro}
 Indeed, if such a resolution
 $$0\rightarrow E^n\rightarrow \ldots
 \rightarrow E^i\rightarrow E^{i-1}\rightarrow E^0\rightarrow {\mathcal I}_Z
 \rightarrow 0$$
 would exist, then we would get the equality
 $$c({\mathcal I}_Z)=\Pi_ic(E_i)^{\epsilon_i}$$
 with $\epsilon_i=(-1)^i$. 
 But the left hand side is non zero in positive degrees
 since its term of degree 
 $r=codim\,Z$ is a non zero multiple of the
 class of $Z$. On the other hand the right hand side
 vanishes in positive degrees by corollary \ref{above}.
 \cqfd
 Note that the assumptions are satisfied
 by a general complex torus  of dimension at least
 $3$. Taking for $Z$ a point, we get  an explicit example
 of a coherent sheaf which does not admit a finite locally free resolution.
 
 \vspace{0,5cm}
 
{\bf Proof of proposition \ref{main}.}
We use again induction on the rank. Let
${\mathcal F}$ be a torsion free coherent sheaf
of rank $k$  on $X$, and
assume first that
${\mathcal F}$ does not contain
any non zero subsheaf of smaller rank.  There is an inclusion
$${\mathcal F}\hookrightarrow {\mathcal F}^{**}$$
whose cokernel is a torsion sheaf, 
where the bidual 
${\mathcal F}^{**}$ of ${\mathcal F}$  is reflexive. Then 
${\mathcal F}^{**}$ does not contain any
non zero subsheaf of smaller rank and hence is stable with respect to the given K\"ahler
metric $h$ on $X$.
The theorem of Bando and Siu together with the
fact that
$$c_1({\mathcal F}^{**})=0,\,<c_2({\mathcal F}^{**}),[\omega]^{n-2}>_X=0$$
  implies that
  ${\mathcal F}^{**}$ is a holomorphic vector bundle
  which is endowed with a flat connection, hence has trivial Chern classes
  and the result is proved in this case.
  
  Assume otherwise that there is an exact sequence
  $$0\rightarrow {\mathcal G}\rightarrow {\mathcal F}\rightarrow {\mathcal H}\rightarrow 0,$$
  where the ranks of
  ${\mathcal G}$ and ${\mathcal H}$ are  smaller than the
  rank of ${\mathcal F}$, and ${\mathcal G}$ and ${\mathcal H}$
  are without torsion.
  This exact sequence determines (and is determined by)
  an extension class
  $$e\in Ext^1({\mathcal H},{\mathcal G}).$$
  Now, by induction on the rank we may assume
  that we have  inclusions
  $${\mathcal G}\hookrightarrow E_1,\,{\mathcal H}\hookrightarrow E_2,$$
  whose cokernels
  ${\mathcal T}_i$  are of torsion and where
  the $E_i$'s are holomorphic vector bundles
  with vanishing Chern classes.
  The extension class
  $e$ gives first an extension class
  $f \in Ext^1({\mathcal H},E_1)$, which provides
  a sheaf
  $E'$ containing
  ${\mathcal F}$ in such way that  $E'/{\mathcal F}$
  is of torsion,
  and fitting into an exact sequence
  $$0\rightarrow E_1\rightarrow E'\rightarrow {\mathcal H}\rightarrow 0.$$
  Next, because  the torsion sheaf ${\mathcal T}_2$ is supported in codimension
  $\geq3$, the restriction map
  provides an isomorphism
  $$Ext^1(E_2,E_1)\cong Ext^1({\mathcal H},E_1)$$
  as one sees by applying Serre's duality to these Ext groups.
  Hence it follows that there is a holomorphic
  vector bundle
  $E$ on $X$, which is an extension of
  $E_2$ by $E_1$, and which contains
  $E'$ as a subsheaf, such that the quotient
  $E/E'$ is of torsion.
  $E$ has vanishing Chern classes, because
  $E_i$ satisfy this property for $i=1,\,2$,
  and contains ${\mathcal F}$ as a subsheaf
  such that the quotient
  $E/{\mathcal F}$ is of torsion. This completes the proof by induction.
  \cqfd

\end{document}